\documentclass[12pt,twoside]{amsart}
\usepackage{amsmath,amssymb,mathrsfs,graphicx,bbm,amsfonts,amsthm}
\usepackage[utf8x]{inputenc}
\usepackage{color}
\usepackage[normalem]{ulem}
\usepackage{todonotes}
\usepackage{comment}

\usepackage{graphicx}

\newtheorem{theorem}{Theorem}[section]
\newtheorem{lemma}[theorem]{Lemma}
\newtheorem{proposition}[theorem]{Proposition}
\newtheorem{corollary}[theorem]{Corollary}
\newtheorem{definition}[theorem]{Definition}

\newtheorem{remark}[theorem]{Remark}

\DeclareMathAlphabet{\mathbfit}{OML}{cmm}{b}{it}

\makeatletter
\@addtoreset{equation}{section}
\makeatother

\newcommand{\cadlag}{\ifmmode\mathcal D\else c\`adl\`ag \fi}

 \newcommand{\ba}{\begin{array}}
 \newcommand{\ea}{\end{array}}
 \newcommand{\bea}{\begin{eqnarray}}
 \newcommand{\eea}{\end{eqnarray}}
 \newcommand{\be}{\begin{equation}}
  \newcommand{\ee}{\end{equation}}

 \def \Z {{\mathbb Z}}
 \def \R {{\mathbb R}}
 
 \def \N {{\mathbb N}}

 \def \E {{\mathbb E}}

 \def \cA {\mathcal{A}}

 \def \cE {\mathcal{E}}
 
 \def \cG {\mathcal{G}}

 \def \cS {\mathcal{S}}

 \def \cW {\mathcal{W}}

 \def \a {{\alpha}}
 \def \b {{\beta}}

 \def \z {{\z}}

 \def \z {{\zeta}}

 \def \L {{\Lambda}}

 \def \1{\mathbbm{1}} % funziona solo con il package bbm.

\begin{document}

\title[]{The GHS and other correlation inequalities\\ for the two-star model}

\author{Alessandra Bianchi}
\address{Dipartimento di Matematica,
Universit\`a degli Studi di Padova, Via Trieste 63,
35121 Padova, Italy.}\email{alessandra.bianchi@unipd.it}
\author{Francesca Collet}
\address{Dipartimento di Informatica,
Universit\`a degli Studi di Verona, Strada le Grazie~15, 
37134 Verona, Italy.}\email{francesca.collet@univr.it}
\author{Elena Magnanini}
\address{WIAS, Mohrenstraße 39, 10117 Berlin.}\email{elena.magnanini@wias-berlin.de}

\date{}

\begin{abstract}
We consider the two-star model, a family of exponential random graphs
indexed by two real parameters, $h$ and $\a$, that rule respectively
the total number of edges and the mutual dependence between them.
Borrowing tools from statistical mechanics, we study different classes of
correlation inequalities for edges, that naturally emerge 
while taking the partial derivatives of the (finite size) free energy. 
In particular, if $\a, h\geq 0$, we derive first and second order 
correlation inequalities and then prove the so-called  GHS inequality.
As a consequence, under the above conditions on the parameters,
the average edge density turns out to be an increasing and concave function of the parameter $h$, at any fixed size of the graph. Some of our results can be extended to more general classes of exponential random graphs.

    \par\smallskip\noindent
    {\it Keywords:} Exponential random graphs; correlation inequalities; Gibbs measure; free energy.
    %\par\smallskip\noindent
    %{\it Acknowledgements:} We want to thank .
\end{abstract}

\maketitle

\section{Introduction}
Correlation inequalities are an important tool in equilibrium statistical mechanics. They are used to estimate moments and correlations in ferromagnetic systems, allowing in turn to obtain analyticity properties of some physical observables (such as magnetization and susceptibility) and to prove/disprove the presence of a phase transition. 
Among these inequalities, we find the Griffiths, Hurst and Sherman (GHS) inequality, that rules the three-particle interactions
and is mainly known for providing convexity properties of relevant functionals.
As the Griffiths, Kelley and Sherman (GKS) inequality \cite{G, GKS},
it was firstly proved for the classical Ising model, to show that the average magnetization 
is a concave function of the positive external field \cite{GHS}, and
 then extended to general classes of even ferromagnets 
that can be derived out of the Ising model \cite{E2, G2, BG, N}.

However, the aforementioned result  is only one of the different implications 
entailed by the GHS inequality. For example, it has been used to characterize possible phase transitions, to prove monotonicity of correlation length, and to derive critical exponent inequalities for the Ising model on $\mathbb{Z}^d$; to obtain monotonicity of mass gap and to estimate coupling constants in $\varphi^4$ field theory; or also to show convexity-preserving properties of certain differential equations and diffusion processes. For further details we refer the reader to \cite{EMN} and references therein.

In the present paper we consider a family of exponential random graphs known 
as two-star model \cite{PN2}. 
Specifically, we consider a Gibbs probability measure on the set of  all simple graphs on 
$n$ vertices, whose Hamiltonian depends on the densities of edges and two-star graphs.
Our goal is to study some correlation inequalities for such a model, 
with a particular focus on the GHS inequality.

In comparison with ferromagnetic systems, the major difference is that
the Gibbs measure of our system, being supported on $\{0,1\}^{\binom{n}{2}}$, 
does not enjoy $\Z_2$-simmetry. 
As a consequence, although the positivity of the support of the measure allows to easily  deduce positivity of the moments and derive
the Fortuin, Kasteleyn and Ginibre (FKG) inequality \cite{FKG},
higher order correlations are non-trivial to analyze, 
and generally depend on the choice of the parameters.

The manuscript is organized as follows.
In Section 2 we introduce the two-star model and we define  
the corresponding free energy function. Moreover, we briefly recall some recent
results about its asymptotic behavior, including the characterization
of the phase diagram and some limit theorems for the edge density.
Section 3 is devoted to correlation inequalities and it collects our 
main results. 
We first provide the formal definition of the aforementioned
FKG, GKS and GHS inequalities in the context of a two-star model
with generalized parameters (see Eq. \eqref{EW_H2}).
In Subsection 3.1 we show that the  FKG and GKS inequalities hold for this model
whenever $\a\geq 0$,
and then we derive some preliminary results used afterwards in the proof
of the GHS inequality, that is the core of the present work (see Theorem \ref{Thm_GHS}).
The statement of this result, that holds under the additional hypothesis $h\geq 0$, is given in Subsection 3.2 together with its proof. 
This is mainly based on ideas from \cite{L}, 
where an alternative and simplified strategy of the original  proof has been devised. 
We then bring back the results to the classical two-star model, and make a few comments about
some immediate consequences of the derived correlation inequalities.
In Section~\ref{Sect:gen} we discuss which of our techniques can be extended to prove the FKG and GKS inequalities for general exponential random graphs and which are the issues in adapting the proofs to obtain the GHS inequality in this setting.

\section{Model and background}
\subsection{Two-star model}
Let us  consider the set $\cG_n$ of all simple graphs on $n$ labeled vertices
that are identified with the elements of the set $[n]=\{1,2,3,\ldots, n\}$.
We define a probability distribution on $\cG_n$
by means of the \textit{homomorphism densities} of the subgraphs of the graph.
Specifically, if $G\in\cG_n$ and $H$ is  a given simple subgraph, we define
\begin{equation}\label{def_graph_hom_density}
t(H,G) := \frac{|\text{hom}(H,G)|}{|V(G)|^{|V(H)|}},
\end{equation}
i.e. the probability that a random mapping $V(H)\mapsto V(G)$
from the vertex set of $H$ to the vertex set of $G$ is edge-preserving.

For any $k \in \mathbb{N}$, let $H_1, H_2, \dots, H_k$ be pre-chosen finite simple graphs (edges, stars, triangles, cycles,~\dots) and let $\boldsymbol{\beta}=(\beta_1, \dots, \beta_k)$ be a collection of real parameters. 
For any choice of $\boldsymbol{\beta}$,
an exponential random graph is identified by the Gibbs probability density
{
\be\label{eq:prob-exprg}
\mu_{n; \boldsymbol{\b}} (G)=\frac{\exp \left(H_{n,\boldsymbol{\b}}(G)\right)} {Z_{n;\boldsymbol{\b}}} \quad \mbox{for } G\in\cG_n.
\ee
}
The function $H_{n,\boldsymbol{\b}}$, called Hamiltonian, is given by
\begin{equation}\label{Hamiltonian}
H_{n;\boldsymbol{\b}}(G)=n^2\sum_{j=1}^{k}\beta_{j}t(H_{j},G)\end{equation}
and the normalizing factor
\be\label{eq:partit-expon}
Z_{n;\boldsymbol{\b}}=\sum_{G \in \cG_{n}} \exp \left(H_{n;\boldsymbol{\b}}(G)\right) 
\ee
is the \textit{partition function}.

In the present setting we focus on the \textit{two-star model},
characterized by a Gibbs measure  that depends only on the densities of edges
and two-star graphs. 
Recall that a two-star graph is an undirected graph with one \textit{root} vertex
and two other vertices connected with the root, and otherwise disconnected.
Under this assumption, the measure can be conveniently expressed as follows.

Let $\mathcal{E}_n$ denote the edge set of the complete graph on $n$ vertices,
with elements  labeled from 1 to $\binom{n}{2}$.
If $i,\,j\in\mathcal E_n$ are neighboring edges, we write $i\sim j$
and we identify the unordered pair $\{i,j\}$ with the resulting two-star graph,
that will be called \textit{wedge} $\{i,j\}$ in short.
Let $\mathcal{W}_n:=\{\{i,j\}: i,j \in \mathcal{E}_n,\, i\sim j\}$ be the set of wedges of $\mathcal E_n$,
and set $\mathcal{A}_n := \{0,1\}^{|\mathcal{E}_n|}$, $|\cdot|$ being the cardinality of a set.

Notice that there is a one-to-one correspondence between graphs
$G\in\cG_n$ and elements ${x}=(x_{i})_{i \in \cE_n}\in\cA_n$ so that,
if $G$ corresponds to $x$, it holds that
%\myinline{Non so se sappiamo il coefficiente del secondo termine}

\be\label{densities}
t(H_1,G)=\frac{2}{n^2} \sum_{i\in\cE_n} x_i \,\qquad
t(H_2,G)=\frac{2}{n^3}\sum_{\{i,j\}\in\cW_n} x_ix_j+\frac{2}{n^3} \sum_{i\in\cE_n} x_i\,,
\ee
with $H_1$ an edge and $H_2$ a wedge.
Hence, we may look at the Hamiltonian of the two-star model as
a function on $\cA_n$  defined by
\begin{equation}\label{EW_H_beta}
H_{n;\beta_1,\beta_2}(x) =
\frac{2\beta_2}{n} \sum_{\{i,j\}\in \mathcal{W}_n} x_i x_j + 2\left(\beta_1+\frac{\beta_2}{n}\right) \sum_{i \in \mathcal{E}_n} x_i\,.
\end{equation}
Notice that this Hamiltonian is asymptotically equivalent (see also \cite{Lo}) to 
\begin{equation}\label{EW_H}
H_{n;\alpha,h}(x) =
\frac{\alpha}{n} \sum_{\{i,j\}\in \mathcal{W}_n} x_i x_j + h \sum_{i \in \mathcal{E}_n} x_i\,,
\end{equation}
where, for convenience, we have set $h=2\b_1$ and $\a=2\b_2$.
In the following, we will focus on the corresponding two-star model, having Gibbs density on $\mathcal{A}_n$
given by
\begin{equation}\label{notaz}
\mu_{n;\alpha,h}(x) =
\frac{\exp\left(H_{n;\alpha,h}(x) \right)}{Z_{n;\alpha,h}} \quad \text{ with }
\quad Z_{n;\alpha,h}= \sum_{x \in \mathcal{A}_n} \exp \left(H_{n;\alpha,h}(x)\right)\,.
\end{equation}
Accordingly, we will denote the related measure and expectation by $\mathbb{P}_{n;\alpha,h}$ and $\E_{n;\alpha,h}$, respectively.

\subsection{Free energy}
The \textit{free energy} is a key function in the context of statistical mechanics, as it  
encodes most of the asymptotic properties of the system. 
Specifically, the finite and infinite size free energies associated with \eqref{EW_H} are
\be\label{free_energy(a,h)}
f_{n;\alpha,h} := \frac{1}{n^2} \ln Z_{n;\alpha,h} \quad \text{ and } \quad
f_{\alpha,h} := \lim_{n\to+\infty}f_{n;\alpha,h}\,.
\ee

To understand the important role of the free energy, we first observe that its partial derivatives w.r.t. $\alpha$ and $h$, respectively give the average edge and wedge densities of the model. 
More precisely, if we denote by $E_n$ the number of edges 
of the graph $G$, and by $W_n$
the number of wedges of $G$,  we get
\be
\label{averages}
\partial_{h} f_{n; \alpha, h} = \frac{\mathbb{E}_{n;\alpha,h}
\left(E_n\right)}{n^2}
\quad \text{ and } \quad
\partial_{\a} f_{n; \alpha,h} = \frac{\mathbb{E}_{n;\alpha,h}
\left(W_n\right)}{n^3}\,.\ee
The characterization of the infinite size free energy, together with its analytical properties, then provides a relevant tool to infer some structural properties of the graph.

As an application of Theorems 4.1 and 6.4 in \cite{CD}, 
for any $(\alpha,h)\in\mathbb{R}^{2}$ we have that 

\begin{equation}\label{free_energy}
f_{\alpha,h} \, = \,  \sup_{0 \leq u \leq 1} \left(\frac{\alpha u^2}{2} + \frac{h u}{2} - \frac{1}{2}I(u) \right) \, = \, \frac{\alpha{(u^{*})}^2}{2} +\frac{h u^{*}}{2} -\frac{1}{2}I(u^{*}),
\end{equation}
where $I(u)= u\ln u + (1-u)\ln (1-u)$ and $u^{*}=u^*(\alpha,h)$ is a maximizer 
that solves the fixed-point equation
\begin{equation}\label{FixPointEq}
\frac{e^{2\alpha\,u +h}}{1+e^{2\alpha\,u +h}}=u.
\end{equation}
Depending on the parameters, Eq \eqref{FixPointEq} can have more than one solution at which the supremum in \eqref{free_energy} is attained. Having multiplicity of optimizers translates in the possibility of having limiting graphs with very different edge densities.

\subsection{Edge-occurrence probability}\label{rem-medie}
As already observed by Park and Newman for the edge-triangle model (see \cite{PN1}, Eq. (4)), 
the probability that the edge $x_i$ is present can be also written as the expectation of a 
function of the Hamiltonian where $x_i=1$.
Explicitly, in our context we obtain
\be\label{PN-gen}
\mathbb{E}_{n;\alpha,h}(x_i)=
\mathbb{E}_{n;\alpha,h} \left[\left(1+\exp\left( 
- \frac{\a}{n}\sum_{j\in\cE_n: j\sim i} x_j -h\right)\right)^{-1}\right]\,.
\ee
Since the model enjoys a symmetry in the edge structure, in the sense that each edge in the complete graph has precisely the same neighborhood, the aforementioned expectation turns out to be the same for all $i$. This leads to
\be\label{banana}
\mathbb{E}_{n;\alpha,h}
\left(E_n\right)= \sum_{i\in\cE_n} \mathbb{E}_{n;\alpha,h}(x_i)= \binom{n}{2}\mathbb{E}_{n;\alpha,h}(x_i)\,.
\ee
Hence, the average edge density  corresponds asymptotically to the edge-occurrence probability. 
\begin{remark}\label{rem-medie3}
At this point, the following remark is in order.
The symmetry in the edge structure is intrinsic to the edge set $\cE_n$,
and does not depend on the specific exponential random graph taken into account.
Hence, the analog of the identity \eqref{banana} holds true for general Hamiltonians 
of the form \eqref{Hamiltonian}. 
To our knowledge, this property, which is evident from the interacting particle system perspective, has not been pointed out before.
\end{remark}

Taking into account identity \eqref{free_energy} and the aforementioned results, 
it holds that
\be
\lim_{n\to\infty}\mathbb{E}_{n;\alpha,h}(x_i) = 
2 \lim_{n\to\infty} \partial_{h} f_{n;\a,h}=
2 \partial_{h} f_{\a,h} = u^*(\a,h)\,.
\ee
While an explicit expression of the edge-occurrence probability
as function of $(\a,h)$ is missing even in the infinite size limit,
it is easy to verify from \eqref{PN-gen} that $\mathbb{E}_{n;\alpha,h}(x_i)\geq 1/2$ for all $n\in\N$, whenever $\a,\,h\geq 0$.
However, simulations suggest that the region of parameters 
where the average edge density is bigger than $1/2$ is larger,  
and it also includes negative values of $h$.
For large enough $n$, this region can be approximately 
characterized by the analysis of the asymptotic behavior of the model.
The study of equations 
\eqref{free_energy} and \eqref{FixPointEq} leads to 
the phase diagram that we are going to summarize.

\subsection{Phase diagram.}  
We collect here the relevant features of the phase diagram of the two-star model, 
that can be obtained as a special case of some of the results in \cite{RY}.
%(see also \cite{AC, AZ, BMPO} for similar statements). 
%
%
 
\begin{figure}[h!]
\centering
\includegraphics[scale=0.6]{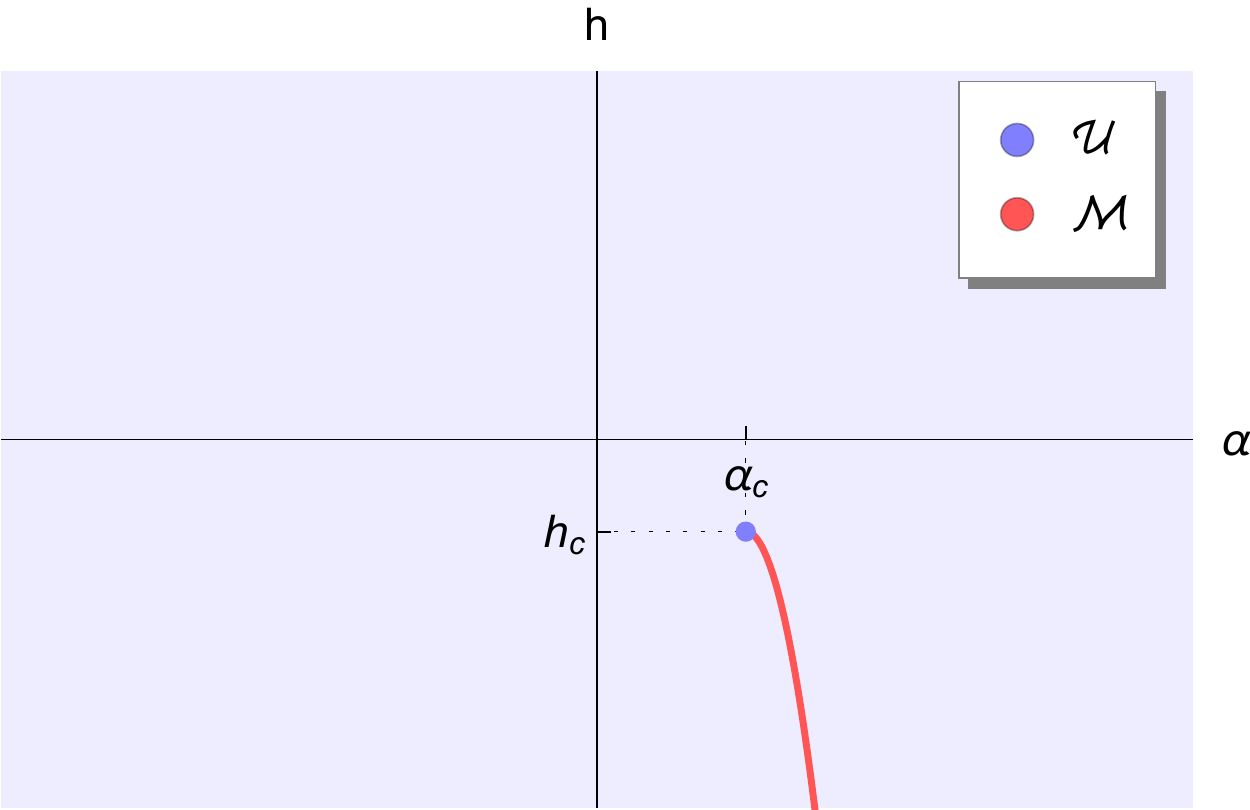}
\caption{\footnotesize Phase space $(\alpha,h)$ for the two-star model \eqref{EW_H}. The blue region, that includes the critical point, is the uniqueness region $\mathcal{U}$ for the maximization problem \eqref{free_energy}; whereas, the red curve corresponds to the critical curve $\mathcal{M}$ along which \eqref{free_energy} admits  two solutions.
}
\label{fig:phase_diagram}
\end{figure}

The infinite size free energy $f_{\alpha,h}$ is well-defined in $\R^2$. 
Moreover, it is analytic in the whole plane except for a continuous critical curve 
\[
%\begin{equation}\label{curve}
\mathcal{M} := \left\{(\alpha,h) \in (\alpha_c,+\infty) \times 
(-\infty,h_{c}): h = q(\alpha)\right\},
%\end{equation}
\]
starting at the critical point $(\alpha_c,h_{c}) = \left(2,-2\right)$ 
and contained in the cone \mbox{$\alpha > 2$}, \mbox{$h < -2$}. 
In particular, the system undergoes a first order phase transition across 
the curve and a second order phase transition at the critical point (see \cite{RY}, Thms.~2.1 \& 2.2). 
The scalar problem \eqref{free_energy} admits one solution in the uniqueness 
region $\mathcal{U} := \mathbb{R}^{2}\setminus \mathcal{M}$ while
it has two solutions along the curve $\mathcal{M}$ (see \cite{RY}, Prop.~3.2).
A qualitative graphical representation of the phase diagram is provided in Fig.~\ref{fig:phase_diagram}.

An analogous analysis has been performed in a sparse regime in \cite{AC}, in the directed graph case in \cite{AZ}, and for a mean-field version of the model in \cite{BMPO}. 
\subsection{Limiting distribution for the edge density}
We summarize some results on the asymptotic behavior of the edge density of the two-star model. By retracing the proofs in \cite{BCM}, we can obtain the following strong law of large numbers and standard central limit theorem: 

\begin{equation}\label{LLN}
\frac{2E_n}{n^2} \, \xrightarrow[n\to\infty]{\;\;\mathrm{a.s.}\;\;}{} \, u^{*}(\alpha,h) \qquad \text{w.r.t. } \mathbb{P}_{n;\alpha,h}, \text{ for }(\alpha,h)\in\mathcal{U}
\end{equation}	
and
\begin{equation}\label{CLT}
\sqrt{2} \, \frac{E_n - \mathbb{E}_{n;\alpha,h}(E_n)}{n} \xrightarrow[n\to\infty]{\;\;\mathrm{d}\;\;} \mathcal{N}(0,v(\alpha,h)) \qquad \text{w.r.t. } \mathbb{P}_{n;\alpha,h}, \text{ for }(\alpha,h)\in\mathcal{U}\setminus\{(\alpha_c,h_c)\},
\end{equation}
where $\mathcal{N}(0,v(\alpha,h))$ is a centered Gaussian distribution with variance 
$v(\alpha,h):=\partial_{h} u^{*}{(\alpha,h)}$, being $u^{*}$ the unique maximizer of \eqref{free_energy}.
%
%The speed of the above convergence has been investigated 
%in \cite{BCM} under different assumption on the parameters.
%
A further result can be also given in the multiplicity region;
for all $(\alpha,h) \in \mathcal{M}$, it holds

\[
\frac{2E_n}{n^2} \, \xrightarrow[n\to\infty]{\;\;\mathrm{d}\;\;}{} \, \kappa \delta_{u_1^{*}(\alpha,h)}+(1-\kappa) \delta_{u_2^{*}(\alpha,h)} \qquad \text{w.r.t. } \mathbb{P}_{n;\alpha,h}, 
\]
where $u_1^*$, $u_2^*$ solve the maximization problem in \eqref{free_energy} and $0 < \kappa < 1$ is a suitable (unknown) constant.

Similar limit theorems are obtained, with different techniques, in \cite{M}, where also results on the partial sum of the degrees can be found.

\section{Correlation inequalities}

In statistical mechanics the study of correlations between particles, so as the analysis of local functions, is often performed with the help of two important inequalities, both related to the sign of the derivatives of the free energy; the \textbf{GKS inequality} and the \textbf{GHS inequality} (see \cite{FV,GHS,GKS,L} and references therein for further details). We aim at deriving the analogs of these two inequalities for
our reference measure $\mu_{n;\a,h}$, given in \eqref{notaz}.

To understand the connection between the GKS inequality and the sign of the
derivatives of the free energy, we introduce  a slightly more general setting. 

Let $\boldsymbol{\alpha} = (\alpha_{ij})_{i,j \in \mathcal{E}_n}$ and
 $\boldsymbol{h}=(h_i)_{i \in \mathcal{E}_n}$ be two collections of real numbers (we write $\boldsymbol{\a} \geq 0$ (resp. $\boldsymbol{h} \geq 0$) as a shortcut for $\a_{ij} \geq 0$ (resp. $h_i \geq 0$) for all $i,j \in \mathcal{E}_n$).
 For $x \in \mathcal{A}_n$, we define the Hamiltonian
\begin{equation}\label{EW_H2}
H_{n;\boldsymbol{\alpha},\boldsymbol{h}}(x) =
\frac{1}{n} \sum_{\{i,j\}\in \mathcal{W}_n}\alpha_{ij} x_i x_j + \sum_{i \in \mathcal{E}_n}h_{i} x_i\,.
\end{equation}
In analogy with \eqref{notaz} and \eqref{free_energy(a,h)}, we denote by $\mu_{n;\boldsymbol{\alpha},\boldsymbol{h}}$ the Gibbs measure obtained from \eqref{EW_H2}, by $\E_{n;\boldsymbol{\alpha},\boldsymbol{h}}$ the corresponding expectation, and we set  $f_{n;\boldsymbol{\alpha},\boldsymbol{h}}
:= \frac{1}{n^2} \ln Z_{n;\boldsymbol{\alpha},\boldsymbol{h}}$ to be the finite size free energy.
Observe that we recover the Hamiltonian \eqref{EW_H} and the related Gibbs measure $\mu_{n;\a,h}$ by setting $\alpha_{ij} \equiv \alpha$, for all $i,j \in \mathcal{E}_n$, and $h_i \equiv h$, for all $i \in \mathcal{E}_n$.

Let $A\subseteq \mathcal{E}_n$ be a given subset of edges. The GKS inequality deals with expectations and covariances of random variables of the type $x_A:=\prod_{i \in A} x_i$, with the convention that $x_{\emptyset}=1$. 

\begin{definition}[GKS inequality]\label{def_GKS}
The Gibbs measure $\mu_{n;\boldsymbol{\a},\boldsymbol{h}}$ on $\mathcal{A}_n$ satisfies the GKS inequality if,
for all $A,B\subseteq \cE_n$,
\be\label{GKS-inequality}
\E_{n;\boldsymbol{\a},\boldsymbol{h}}(x_A x_B)
\geq \E_{n;\boldsymbol{\a},\boldsymbol{h}}(x_A)\cdot
\E_{n;\boldsymbol{\a},\boldsymbol{h}}(x_B)\,.
\ee
\end{definition}
\begin{remark}
Notice that, by choosing  $A=\{i\}$ and $B=\{j\}$, with $i\neq j$,
from the GKS inequality it follows  that $x_i$ and $x_j$
are positively  correlated under $\mu_{n;\boldsymbol{\a},\boldsymbol{h}}$.
\end{remark}

A useful link between the correlations of the system and the partial derivatives of the free energy w.r.t. the parameters $h_i$'s is provided by the MacLaurin expansion of the 
log moment generating function of $x \in \mathcal{A}_n$.

The coefficients of this expansion are the so-called \emph{Ursell functions},
that are formally defined,  for $\ell \in [n]$ and  any choice of
 $i_1, \dots, i_{\ell} \in \mathcal{E}_n$, by
\be\label{ursell}
u_{\ell}(i_1, \dots,i_{\ell}) := n^2\partial_{h_{i_1}\ldots h_{i_\ell}}
f_{n;\boldsymbol{\alpha},\boldsymbol{h}}\,.
\ee
For instance, this yields
\begin{align}
u_1(i) &= \E_{n;\boldsymbol{\alpha},\boldsymbol{h}} (x_i ),\\
u_2(i,j) &= \E_{n;\boldsymbol{\alpha},\boldsymbol{h}}( x_i x_j ) -
\E_{n;\boldsymbol{\alpha},\boldsymbol{h}}( x_i )\E_{n;\boldsymbol{\alpha},\boldsymbol{h}}( x_j ),\\
u_3(i,j,k) &= \E_{n;\boldsymbol{\alpha},\boldsymbol{h}}(  x_i x_j x_k ) -
\E_{n;\boldsymbol{\alpha},\boldsymbol{h}}( x_i ) \, \E_{n;\boldsymbol{\alpha},\boldsymbol{h}}( x_j x_k )
-
\E_{n;\boldsymbol{\alpha},\boldsymbol{h}}( x_j ) \,
\E_{n;\boldsymbol{\alpha},\boldsymbol{h}}( x_i x_k )
 \notag\\[.2cm]
 \label{u3}
&- \E_{n;\boldsymbol{\alpha},\boldsymbol{h}}( x_k ) \,
\E_{n;\boldsymbol{\alpha},\boldsymbol{h}}( x_i x_j )
 + 2 \E_{n;\boldsymbol{\alpha},\boldsymbol{h}}( x_i )
 \E_{n;\boldsymbol{\alpha},\boldsymbol{h}}( x_j )
 \E_{n;\boldsymbol{\alpha},\boldsymbol{h}}( x_k ).
\end{align}

\begin{remark}\label{rem-Ursell}
Notice that the definition of the Ursell functions \eqref{ursell} 
necessarily passes through the generalized  setting with vector parameters 
$\boldsymbol{\alpha},\boldsymbol{h}$, of which they are functions.
However, when computed along the constant vectors $\boldsymbol{\alpha}$ and $\boldsymbol{h}$,
with $\a_{ij}\equiv \a$ for all $i,j \in \mathcal{E}_n$, and $h_i\equiv h$ 
for all $i \in \mathcal{E}_n$, they are also useful to characterize  the derivatives 
of the classical free energy $f_{n;\a,h}$ through the identities
\be\label{uresell2}
n^2 \partial_{\underbrace{\text{\scriptsize $h\ldots h$}}_{\text{$\ell$ times}}} f_{n;\a,h}= \sum_{i_1, \dots,i_{\ell}\in\cE_n} u_{\ell}(i_1, \dots,i_{\ell})\,,\qquad \forall \ell \in [n]\,.
\ee
%$$
% \partial_{\underbrace{\text{\scriptsize $h\ldots h$}}_{\text{$\ell$ times}}} 
%f_{n;\a,h}= \sum_{i_1, \dots,i_{\ell}\in\cE_n}
%\partial_{h_{i_1}\ldots h_{i_\ell}}
%f_{n;\boldsymbol{\alpha},\boldsymbol{h}}|_{\boldsymbol{\alpha}=\a;\boldsymbol{h}=h}\,.
%$$

\end{remark}

While the GKS inequality implies \mbox{$u_2(i,j)\geq 0$}, giving positive correlation between the random variables $x_i$ and $x_j$, the GHS inequality concerns  the sign of the Ursell function $u_3(i,j,k)$.

\begin{definition}[GHS inequality]\label{def_GHS}
The Gibbs measure $\mu_{n;\boldsymbol{\alpha},\boldsymbol{h}}$ on $\mathcal{A}_n$
satisfies the GHS inequality if,
for all $i,j,k\in\cE_n$, $u_3(i,j,k)\leq 0$
or, equivalently, if
\begin{multline}\label{GHS-inequality}
\E_{n;\boldsymbol{\alpha},\boldsymbol{h}}(  x_i x_j x_k ) -
\E_{n;\boldsymbol{\alpha},\boldsymbol{h}}( x_i ) \, 
\E_{n;\boldsymbol{\alpha},\boldsymbol{h}}( x_j x_k )
-
\E_{n;\boldsymbol{\alpha},\boldsymbol{h}}( x_j ) \,
\E_{n;\boldsymbol{\alpha},\boldsymbol{h}}( x_i x_k )\\
- \E_{n;\boldsymbol{\alpha},\boldsymbol{h}}( x_k ) \,
\E_{n;\boldsymbol{\alpha},\boldsymbol{h}}( x_i x_j )
 + 2 \E_{n;\boldsymbol{\alpha},\boldsymbol{h}}( x_i )
 \E_{n;\boldsymbol{\alpha},\boldsymbol{h}}( x_j )
 \E_{n;\boldsymbol{\alpha},\boldsymbol{h}}( x_k )
 \leq 0.
\end{multline}
\end{definition}

Observe that, in our case, $u_1(i) \geq 0$ trivially, due to the fact that $x_i \in \{0,1\}$. The rest of the section is devoted to proving $u_2(i,j) \geq 0$ and $u_3(i,j,k) \leq 0$.

\subsection{The FKG and GKS inequalities}

We start with a preliminary result, the FKG inequality, that will help us in deriving the more advanced inequalities \eqref{GKS-inequality} and~\eqref{GHS-inequality}.

We first show that the measure $\mu_{n;\boldsymbol{\alpha},\boldsymbol{h}}$ on $\mathcal{A}_n$ satisfies a proper lattice condition. Recall that $\mathcal{A}_n$ is partially ordered by 
\begin{equation}\label{partial_order_definition}
x \leq y \quad \text{ if } \quad x_i \leq y_i \quad \text{ for all } i \in \mathcal{E}_n.
\end{equation}
Moreover, given two configurations $x, y \in \mathcal{A}_n$, the (pointwise) maximum and minimum configurations are defined as
\[
(x \vee y)(i) := \max \{x_i, y_i\} \quad \text{ and } \quad (x \wedge y)(i) := \min \{x_i, y_i\},
\]
for all $i \in \mathcal{E}_n$.
The following property holds true.

\begin{lemma}\label{lat_cond}
If $\boldsymbol{\a}\geq 0$, then the Gibbs measure
$\mu_{n;\boldsymbol{\alpha},\boldsymbol{h}}$
fulfills the FKG lattice condition
\begin{equation}\label{FKG_lattice_condition_ETmodel}
\mu_{n;\boldsymbol{\alpha},\boldsymbol{h}} (x \vee y) \, \mu_{n;\boldsymbol{\alpha},\boldsymbol{h}} (x \wedge y) \geq \mu_{n;\boldsymbol{\alpha},\boldsymbol{h}} (x) \, \mu_{n;\boldsymbol{\alpha},\boldsymbol{h}} (y)
\quad \text{ for } x,y \in \mathcal{A}_n.
\end{equation}
\end{lemma}

\begin{proof}
For a configuration $z \in \mathcal{A}_n$, let $E_{z}:=\{i\in\cE_n:\, z_i=1\}$, namely the set of edges present in $z$.
 We have $E_{x \vee y} = E_x \cup E_y$ and $E_{x \wedge y} = E_x \cap E_y$. Observe that
 \begin{itemize}
\item the edges in the configuration $x \vee y $ are the edges the configurations $x$ and $y$ have in common, the edges present in configuration $x$ only and those present in configuration $y$ only;
\item the edges in the configuration $x \wedge y $ are the edges the configurations $x$ and $y$ have in common;
\item the wedges in the configuration $x \vee y $ are the wedges the configurations $x$ and $y$ have in common, the wedges present in configuration $x$ (resp. configuration $y$) only and the wedges you may create by superimposing the edges of the two configurations;
\item the wedges in the configuration $x \wedge y$ are the wedges the configurations $x$ and $y$ have in common.
\end{itemize}
Therefore, verifying that \eqref{FKG_lattice_condition_ETmodel}
is satisfied reduces to show the validity of the inequality
\begin{equation}\label{cc-FKG}
\exp \left\{ \frac{1}{n} \sum_{\{i,j\} \in E} \alpha_{ij} x_i x_j \right\} \geq 1,
\end{equation}
where
\[
E = \left\{ \{i,j\}: \{i, j\} \subset E_{x \vee y} \text{ is a wedge and } \{i,j\} \left[ \begin{scriptsize} \begin{array}{l} \not\subset E_x \\ \not\subset E_y \end{array} \end{scriptsize} \right.\right\}.
\]
The conclusion follows as $\boldsymbol{\alpha} \geq 0$ by assumption.
\end{proof}

An immediate consequence of Lemma \ref{lat_cond} is the positive correlation of increasing random variables. Specifically,  if $f$ and $g$  are increasing functions on $\cA_n$ (i.e., \mbox{$f(x)\leq f(y)$} if $x\leq y$), then we obtain
the \textbf{FKG inequality}

\be\label{FKG inequality}
\E_{n;\boldsymbol{\alpha},\boldsymbol{h}}(f g)
\geq \E_{n;\boldsymbol{\alpha},\boldsymbol{h}}(f)\cdot
\E_{n;\boldsymbol{\alpha},\boldsymbol{h}}(g)\,.
\ee

\begin{corollary}
If $\boldsymbol{\a}\geq 0$, then the Gibbs measure
$\mu_{n;\boldsymbol{\alpha},\boldsymbol{h}}$ 
satisfies the GKS inequality.
\end{corollary}

\begin{proof}
Notice that for all $A\subseteq \mathcal{E}_n$,
the function $x_A=\prod_{i \in A} x_i$ is increasing in \mbox{$x \in \mathcal{A}_n$}.
Hence, by applying the FKG inequality \eqref{FKG inequality}
to the functions $f(x)=x_A$ and $g(x)=x_B$,
we immediately derive \eqref{GKS-inequality}.
\end{proof}

\begin{remark}\label{FKG-GKS-wedgemodel}
A straightforward adaptation of the arguments of Lemma \ref{lat_cond}
applies to general exponential random graphs. 
We refer the reader to Section \ref{Sect:gen} for more details.
\end{remark}

We now provide two useful consequences of the GKS inequality.
To state properly the results we need to introduce a few more notation; we need a suitable ``restriction'' of the system
%measure $\mu_{n; \boldsymbol{\alpha},\boldsymbol{h}}$ 
on a subset.

For $A\subseteq\cE_n$, set $\mathcal{W}_A:=\{\{i,j\}: i,j \in A\,,\, i\sim j\}$
and define the Hamiltonian 
 \be\label{ham-subset}
H_{A;\boldsymbol{\alpha},\boldsymbol{h}}(x)
= \frac{1}{n} \sum_{\{i,j\} \in \mathcal{W}_A} \alpha_{ij} x_i x_j
+ \sum_{i \in A} h_i x_i\,.
\ee
Let $\mu_{A; \boldsymbol{\alpha},\boldsymbol{h}}$ be the associated Gibbs measure, with normalizing constant
$Z_{A; \boldsymbol{\alpha},\boldsymbol{h}}$ (partition function), and let $\E_{A; \boldsymbol{\alpha},\boldsymbol{h}}$ denote the corresponding  expectation.

A first consequence of the GKS inequality is a form of monotonicity, with respect to the volume, that can be established for the averages of  $x_\L$, with $\L\subseteq \cE_n$.

\begin{lemma}\label{lem_GKS-volume}
If the Gibbs measure $\mu_{n;\boldsymbol{\alpha},\boldsymbol{h}}$ satisfies the GKS inequality then, for any $\L\subseteq A\subseteq B\subseteq \cE_n$,
\be\label{GKS-wedge2}
\E_{A; \boldsymbol{\alpha},\boldsymbol{h}}( x_\L)\leq
\E_{B; \boldsymbol{\alpha},\boldsymbol{h}}(x_\L)\, .
\ee
\end{lemma}

\begin{proof}
Observe first that, for all $\L\subseteq A\subseteq \cE_n$,
the function $\E_{A; \boldsymbol{\alpha},\boldsymbol{h}} (x_\L)$ is non-decreasing
in $\boldsymbol{\a}$.
Indeed, by differentiating  $\E_{A;\boldsymbol{\alpha},\boldsymbol{h}}$
w.r.t $\a_{ij}$, we get
 \be\label{crescenza-param}
\partial_{\a_{ij}}
 \E_{A;\boldsymbol{\alpha},\boldsymbol{h}}(x_\L)
=\frac{1}{n}\left(\E_{A;\boldsymbol{\alpha},\boldsymbol{h}}(x_\L x_i x_j)-
\E_{A;\boldsymbol{\alpha},\boldsymbol{h}}(x_\L)
\E_{A;\boldsymbol{\alpha},\boldsymbol{h}}(x_i x_j)\right)
\geq 0 \,,
\ee
where the last inequality follows  from the GKS inequality.

Now let $\cW_{A,B}:=\{ \{i,j\}: i\in A \,, j\in B \setminus A\,,\, i\sim j\}$
and, for $s\in[0,1]$, consider the Hamiltonian
$$
H_{B;\boldsymbol{\alpha}(s),\boldsymbol{h}}(x):= \frac{1}{n}\sum_{\{i,j\}\in \cW_B\setminus \cW_{A,B}} \alpha_{ij} x_ix_j
+\frac{s}{n} \sum_{\{i,j\}\in\cW_{A,B}} \alpha_{ij} x_ix_j + \sum_{i\in B} h_ix_i\,,
$$
with corresponding Gibbs measure $\mu_{B;\boldsymbol{\alpha}(s),\boldsymbol{h}}$ and relative average
$\E_{B;\boldsymbol{\alpha}(s),\boldsymbol{h}}$.
Notice that, if $s=1$, we obtain the system on the set $B$, so that $\E_{B;\boldsymbol{\alpha}(1),\boldsymbol{h}}( x_\L)=
\E_{B; \boldsymbol{\alpha},\boldsymbol{h}}( x_\L)$. 
Moreover, since  $\cW_B= \cW_A \sqcup \cW_{B\setminus A} \sqcup \cW_{A,B}$,
when $s=0$, we get
$$
H_{B;\boldsymbol{\alpha}(0),\boldsymbol{h}}(x)
= H_{A;\boldsymbol{\alpha},\boldsymbol{h}}(x) +
H_{B\setminus A;\boldsymbol{\alpha},\boldsymbol{h}}(x)\,.
$$
Then $\mu_{B;\boldsymbol{\alpha}(0),\boldsymbol{h}}=
\mu_{A;\boldsymbol{\alpha},\boldsymbol{h}}
\cdot \mu_{B\setminus A;\boldsymbol{\alpha},\boldsymbol{h}}$
and, as a consequence, being $\L\subseteq A$, we have
$
\E_{B;\boldsymbol{\alpha}(0),\boldsymbol{h}}( x_\L)
=\E_{A; \boldsymbol{\alpha},\boldsymbol{h}}( x_\L)$.
%$\qquad
%\E_{B;\boldsymbol{\alpha},\boldsymbol{h}}^{1}( x_\L)=
%\E_{B; \boldsymbol{\alpha},\boldsymbol{h}}( x_\L)\,.
%$
Finally, since $\boldsymbol{\a} \mapsto \E_{B; \boldsymbol{\alpha},\boldsymbol{h}}(x_\L)$ is 
a non-decreasing mapping and $\boldsymbol{\a}(0)<\boldsymbol{\a}(1)$,
we conclude
$$
\E_{A; \boldsymbol{\alpha},\boldsymbol{h}}( x_\L)=
\E_{B; \boldsymbol{\alpha}(0),\boldsymbol{h}}( x_\L)
\leq \E_{B; \boldsymbol{\alpha}(1),\boldsymbol{h}}( x_\L)
= \E_{B; \boldsymbol{\alpha},\boldsymbol{h}}( x_\L)\,,$$
as claimed.
\end{proof}

A second consequence of the GKS inequality is a comparison between
partition functions.

\begin{lemma}\label{lem_zeta}
If the Gibbs measure $\mu_{n;\boldsymbol{\alpha},\boldsymbol{h}}$ satisfies the GKS inequality then, for any
$E,F \subseteq \cE_n$,
\be\label{zeta}
Z_{E;\boldsymbol{\alpha},\boldsymbol{h}}Z_{F;\boldsymbol{\alpha},\boldsymbol{h}}
\leq
Z_{E\cup F;\boldsymbol{\alpha},\boldsymbol{h}}Z_{E\cap F;\boldsymbol{\alpha},\boldsymbol{h}}\,.
\ee
\end{lemma}
\begin{proof}
We follow some ideas developed in \cite{L} to prove
an analogous result for Ising spin systems. We set
$K_1:=E\cap F$, $K_2:= E\setminus K_1$ and $K_3:= F\setminus K_1$,
so that we can express the sets $E$, $F$, $E\cup F$ and $E\cap F$
as  proper disjoint unions of the subsets $K_i$'s.
With this notation, the inequality \eqref{zeta} becomes equivalent to
 \be\label{prova1}
 L(\boldsymbol{\alpha},\boldsymbol{h}):=\ln Z_{K_1\cup K_2\cup K_3;\boldsymbol{\alpha},\boldsymbol{h}}
-\ln \frac{Z_{K_1\cup K_2;\boldsymbol{\alpha},\boldsymbol{h}}Z_{K_1\cup K_3;\boldsymbol{\alpha},\boldsymbol{h}}}
{Z_{K_1;\boldsymbol{\alpha},\boldsymbol{h}}}
\geq 0\,.
\ee

Notice that, if there is no interaction between the edges in $K_1$ and those in $K_3$,
then
$$Z_{K_1\cup K_2\cup K_3;\boldsymbol{\alpha},\boldsymbol{h}}
=Z_{K_1\cup K_2;\boldsymbol{\alpha},\boldsymbol{h}}Z_{K_3;\boldsymbol{\alpha},\boldsymbol{h}} \quad \text{ and } \quad
Z_{K_1\cup K_3;\boldsymbol{\alpha},\boldsymbol{h}}=Z_{K_1;\boldsymbol{\alpha},\boldsymbol{h}}Z_{K_3;\boldsymbol{\alpha},\boldsymbol{h}},$$
that yields $L(\boldsymbol{\alpha},\boldsymbol{h})=0$. To conclude, it suffices to show that the function $L(\boldsymbol{\alpha},\boldsymbol{h})$
is non-decreasing with respect to the interaction parameter $\boldsymbol{\a}$.
To this purpose, we consider the change  in $L(\boldsymbol{\alpha},\boldsymbol{h})$,
when an interaction of strength $\a_{ij}$, between the edges $i\in K_1$ and $j\in K_3$, 
is added to the system. 
By differentiating w.r.t. $\a_{ij}$ we get
\be\label{prova2}
\partial_{\a_{ij}}
L(\boldsymbol{\a},\boldsymbol{h}) =
\frac{1}{n}\left(\E_{K_1\cup K_2 \cup K_3; \boldsymbol{\a}, \boldsymbol{h}}(x_ix_j)
-\E_{K_1\cup K_3; \boldsymbol{\a}, \boldsymbol{h}}(x_ix_j)\right)\geq 0\,,
\ee
where the last inequality follows from Lemma \ref{lem_GKS-volume}.
All together this implies that \mbox{$L(\boldsymbol{\alpha},\boldsymbol{h})\geq 0$}.
\end{proof}

%\begin{remark}\label{FKG-GKS-wedgemodel}
%A straightforward adaptation of the arguments of Lemma \ref{lat_cond}
% applies to the edge-triangle model with Hamiltonian
%%
%\begin{equation}\label{Hamilt_ERG}
%H_{n;\alpha,h}(x) = \frac{\alpha}{n} \sum_{\{i,j,k\} \in \mathcal{T}_n} x_i x_j x_k + h \sum_{i \in \mathcal{E}_n} x_i, \quad \alpha\geq\,0, h\in\mathbb{R}
%\end{equation}
%%
%where $\mathcal{T}_n=\{\{i,j,k\} \subset \mathcal{E}_n: \{i,j,k\} \text{ is a triangle}\}$.
%This implies that whenever $\a\geq 0$, and for all $h\in\R$, the
% FKG and GKS inequalities w.r.t. the Gibbs measure associated to \eqref{Hamilt_ERG}
% are verified.
%\end{remark}

\subsection{The GHS inequality}\label{subsec-GHS}
We are now ready to derive our main result: the GHS inequality for the model associated with the Hamiltonian \eqref{EW_H2}.

\begin{theorem}\label{Thm_GHS}
If $\boldsymbol{\a},\boldsymbol{h}\geq 0$,
 then the Gibbs measure
$\mu_{n;\boldsymbol{\a},\boldsymbol{h}}$ satisfies the GHS inequality. In particular, for any choice of $i,j,k\in\mathcal{E}_n$, we have
\be\label{GHS-generale}
\partial_{h_{i} h_{j} h_{k}} f_{n;\boldsymbol{\a},\boldsymbol{h}}\leq 0\,.
\ee
\end{theorem}
\begin{remark}
The above theorem  provides sufficient conditions for the validity
of the GHS inequality, and it is then natural to wonder whether
they are also necessary.
A hint on this question is given when considering the statement 
for indices $i=j=k$. Inequality \eqref{GHS-generale} reduces to  
(see also \eqref{GHS-inequality}) 
\be\label{GHS-i=j=k}
\E_{n;\boldsymbol{\alpha},\boldsymbol{h}}(x_i)
\left(1- \E_{n;\boldsymbol{\alpha},\boldsymbol{h}}(x_i)\right)
\left(1- 2\E_{n;\boldsymbol{\alpha},\boldsymbol{h}}(x_i)\right)
\leq 0\,,
\ee
that is verified if and only if
$$\E_{n;\boldsymbol{\alpha},\boldsymbol{h}}(x_i)\geq 1/2\,.$$
Recall that by \eqref{PN-gen}  the above condition is fulfilled whenever 
$\boldsymbol{\a},\boldsymbol{h}\geq 0$. 
This assumption is indeed the only strict requirement on the parameter 
$\boldsymbol{h}$ that we will use along the proof, and precisely in \eqref{passo2} below,
though in a modified setting that requires the validity of this condition uniformly in $n$.
However, as mentioned in Subsection \ref{rem-medie}, the edge-occurence probabilities 
are implicit functions of the parameters $\boldsymbol{\a}$
and $\boldsymbol{h}$, and are also dependent on $n$.  
For these reasons, we believe that the derivation of explicit necessary conditions 
could be in general a hard task.

%To state the following result, it is convenient introduce the  
%set of parameters 
%\be
%\cR^+:=\{ (\boldsymbol{\a},\boldsymbol{h})\,:\, \boldsymbol{\a}\geq 0
%\mbox{ and }
%\E_{n;\boldsymbol{\alpha},\boldsymbol{h}}(x_i)\geq \frac 12 \,,\, \forall i\in\cE_n
%  \}
%\,.\ee
\end{remark}

The strategy of the proof is based on the trick of introducing a duplicate set of variables.
Let $y \in \mathcal{A}_n$ be an independent copy of
$x \in \mathcal{A}_n$, with the same Hamiltonian as in \eqref{EW_H2}, and let $\E$ denote the expectation with respect to the joint measure
\begin{equation}\label{ET_joint_measure}
\mu(x,y) :=
\frac{\exp \left\{H_{n;\boldsymbol{\a},\boldsymbol{h}}(x) +H_{n;\boldsymbol{\a},\boldsymbol{h}}(y) \right\}}{Z_{n;\boldsymbol{\a},\boldsymbol{h}}^2}.
\end{equation}
For any $i \in \mathcal{E}_n$, define the variables $z_i = x_i - y_i$ and
$v_i = \frac{1}{2}(x_i + y_i)$.
Notice that $z_i$ takes value on $\{-1,0,+1\}$, while $v_i$ takes value on
$\left\{0,\frac{1}{2},1\right\}$, and that the following equivalence of events holds
for all $i \in \mathcal{E}_n$:
\be\label{eventi}
\left\{z_i \in \{-1,+1\}\right\} = \left\{v_i = \frac{1}{2}\right\}\quad \text{ and } \quad
\left\{v_i  \in \{0,1\}\right\}=\left\{ z_i =0\right\}.
\ee
With standard notation we set  $z:=(z_i)_{i\in\cE_n}$ and $v:=(v_i)_{i\in\cE_n}$. Moreover, for any given $A\subseteq \mathcal{E}_n$, we define the functions $z_A:=\prod_{i \in A} z_i$ and $v_A:=\prod_{i \in A}v_i$.

\begin{proposition}\label{Prop:GHS-zv}
Let $\boldsymbol{\alpha}, \boldsymbol{h}\geq 0$.
Then, for any $C,D \subseteq \mathcal{E}_n$, it holds that
\begin{align}
\E( z_C z_D)&\geq \E( z_C)\E( z_D )\,, \label{corr_in_2EG}
\\[.2cm]
\E(z_C v_D )&\leq \E(z_C )\E( v_D )\,.
\label{corr_in_1EG}
\end{align}
\end{proposition}

\begin{remark}\label{chiusura}
It is easy to check that the  Ursell function $u_3(i,j,k)$, given explicitly
in \eqref{u3},  can be written as a function of the random variables $z_i$'s and $v_i$'s as
\begin{equation}\label{GHS_inequality_v2}
u_3(i,j,k) = \E( z_i z_j v_k)-\E(z_i z_j)
\E( v_k )\,.
\end{equation}
The statement of Theorem \ref{Thm_GHS} is then a consequence
of the inequality \eqref{corr_in_1EG}.
Similarly, it can be shown that Eq. \eqref{corr_in_2EG} implies
the GKS inequality for the Gibbs measure $\mu_{n;\boldsymbol{\a},\boldsymbol{h}}$.
\end{remark}

%\subsubsection{Proof of the GHS inequality}
\begin{proof}[Proof of Proposition \ref{Prop:GHS-zv}]
We first consider two general functions $\Phi(z)$ and $\Psi(v)$,
with $z=(z_i)_{i\in \mathcal E_n}$  and $v=(v_i)_{i\in \mathcal E_n}$, and we try to express the average $\E(\Phi(z)\Psi(v))$ in a convenient form. Later
we will focus on the functions $\Phi(z)= z_C$ and $\Psi(v)=v_D$.

Observe that, due to the identity $x_ix_j + y_iy_j = \frac 12 z_iz_j + 2 v_iv_j$,
the exponent of the joint measure \eqref{ET_joint_measure}
can be phrased in terms of the variables $z$ and $v$. It yields
\be\label{splitting}
H_{n;\boldsymbol{\alpha},\boldsymbol{h}}(x) +H_{n;\boldsymbol{\alpha},\boldsymbol{h}}(y) = \widehat{H^1}_{n;\boldsymbol{\alpha}}(z)+
\widehat{H^2}_{n;\boldsymbol{\alpha},\boldsymbol{h}}(v),
\ee
where
\begin{equation}\label{Hamil-zv}
\begin{split}
&\widehat{H^1}_{n;\boldsymbol{\alpha}}(z) = \frac{1}{2n}
\sum_{\{i,j\} \in \mathcal{W}_n} \alpha_{ij} z_i z_j\,,\\
&\widehat{H^2}_{n;\boldsymbol{\alpha},\boldsymbol{h}}(v) = \frac{2}{n}
\sum_{\{i,j\} \in \mathcal{W}_n} \alpha_{ij}v_i v_j +
2 \sum_{i \in \mathcal{E}_n} h_i v_i\,.
\end{split}
\end{equation}
Moreover, by exploiting the constraints \eqref{eventi}, we can partition the state space of the couple $(z,v)$ in a disjoint union, over subsets
$A\subseteq \mathcal E_n$, of the sets
\be\label{partizione}
\cS_A := \left\{(z,v): z_i=0, v_i\in \{0,1\} \, \forall i\in A \mbox{ and }
v_i=\frac{1}{2}, z_i\in \{-1,1\}  \, \forall i\in A^c\right\}.
\ee
%$$\cS_A :=\{(z,v)\,:\,z_i=0 \Longleftrightarrow i\in A \}
%=\{(z,v)\in \cA_n^2\,:\,v_i=1/2 \Longleftrightarrow i\in A^c \}\,.$$
Hence, we can write
\begin{equation}\label{mediafunzioni}
\E(\Phi(z)\Psi(v))=
\sum_{A\subseteq \mathcal E_n}\sum_{(z,v)\in \mathcal S_A}
\Phi(z)\Psi(v)
\frac{\exp\left\{ \widehat{H^1}_{n;\boldsymbol{\alpha}}(z)+\widehat{H^2}_{n;\boldsymbol{\alpha},\boldsymbol{h}}(v)\right\}}{Z_{n;\boldsymbol{\alpha},\boldsymbol{h}}^2}\,.
\end{equation}

It is easy to see that if $(z,v)\in\mathcal S_A$, and with the same notation
introduced in \eqref{ham-subset}, we obtain

\begin{equation}\label{hamil-belle-Ising}
\widehat{H^1}_{n;\boldsymbol{\alpha}}(z)
= \frac{1}{2n}
\sum_{ \{i,j\} \in \mathcal{W}_{A^c}} \alpha_{ij} z_i z_j\,, \quad \mbox{with }
z_i\in\{-1,1\}, \,\forall i\in A^c
\end{equation}
and
\begin{multline}\label{hamil-belle-star}
\widehat{H^2}_{n;\boldsymbol{\alpha},\boldsymbol{h}}(v)
= \frac{2}{n}
\sum_{\{i,j\} \in \mathcal{W}_A} \alpha_{ij} v_i v_j + \sum_{i \in A} \left( 2h_i + \frac{1}{n} \sum_{j \in A^c: j \sim i} \alpha_{ij} \right) v_i \\
+ \frac{1}{2n} \sum_{ \{i,j\} \in \mathcal{W}_{A^c}} \alpha_{ij} + \sum_{i \in A^c} h_i
\,, \qquad \mbox{with }
v_i\in\{0,1\}, \,\forall i\in A\,.
\end{multline}
In particular, on the set $\mathcal S_A$, 

\begin{itemize}
\item
the Hamiltonian $\widehat{H^1}_{n;\boldsymbol{\alpha}}(z)$
corresponds to the Hamiltonian of an Ising spin system on the set $A^c$,
with inverse temperature $\boldsymbol{\beta} := \boldsymbol{\alpha}/2n$, 
magnetic field $\boldsymbol{h}=\boldsymbol{0}$, and associated Gibbs
measure
$$\mu_{A^c;\boldsymbol{\beta},\boldsymbol{0}}^{\mathrm{Is}}(z):= \frac{e^{H_{A^c;\boldsymbol{\b},\boldsymbol{0}}^{\mathrm{Is}}(z)}}
{Z_{A^c;\boldsymbol{\b},\boldsymbol{0}}^{\mathrm{Is}}}\,.$$
\item
the Hamiltonian $\widehat{H^2}_{n;\boldsymbol{\alpha},\boldsymbol{h}}(v)$
corresponds to the two-star Hamiltonian on $A$ given in \eqref{ham-subset}, but with parameters
$\boldsymbol{\alpha'} := 2\boldsymbol{\alpha}$
and $\boldsymbol{h'} := (h_i')_{i \in \mathcal{E}_n}$, where $h_i' := 2h_i + \frac{1}{n} \sum_{j \in A^c: j \sim i} \alpha_{ij}$. Indeed, the two Hamiltonians differ only for the constant term $\frac{1}{2n} \sum_{\{i,j\} \in \mathcal{W}_{A^c}} \alpha_{ij} + \sum_{i \in A^c} h_i$
that, being irrelevant in the Gibbs measure, will be neglected. As before, we write $\mu_{A;\boldsymbol{\alpha'},\boldsymbol{h'}}$ for the Gibbs measure related to the Hamiltonian \eqref{hamil-belle-star}.
 \end{itemize}

Going back to Eq. \eqref{mediafunzioni}, in view of the previous considerations, it turns out that
\begin{equation}\label{cambiospazio}
\E( \Phi(z)\Psi(v))=
\sum_{A\subseteq \mathcal E_n}
P(A) f^{\Phi}(A)g^{\Psi}(A),
\end{equation}
where, with self-explanatory notation, we set
\be
f^{\Phi}(A):= \E_{A^c, \boldsymbol{\b},\boldsymbol{0}}^{\mathrm{Is}}(\Phi(z)|_{z_i=0, \forall i \in A})\,,\qquad
g^{\Psi}(A):= \E_{A,\boldsymbol{\alpha'},\boldsymbol{h'}}(\Psi(v)|_{v_{i}=\frac 12 , \forall i \in A^c})
\ee
and
\be\label{misura_P}
P(A):=\frac{Z_{A^c;\boldsymbol{\beta},\boldsymbol{0}}^{\mathrm{Is}}\cdot
Z_{A;\boldsymbol{\alpha'},\boldsymbol{h'}}}{Z_{n;\boldsymbol{\alpha},\boldsymbol{h}}^2}\,.
\ee
Notice that $P$ is a probability on $\cE_n$ by construction.
%Indeed, by \eqref{partizione},
%$$
%\sum_{A\subseteq \cE_n} P(A)=
%\sum_{A\subseteq \cE_n} \sum_{(z,v)\in\cS_A} \frac{\exp\left\{ \widehat{H^1}_{n;\alpha}(z)+\widehat{H^2}_{n;\alpha,h}(v)\right\}}{\left(Z_{n;\alpha,h}\right)^2}=1\,.
%$$
Specializing  the identity \eqref{cambiospazio} to the functions
$\Phi(z)=z_C$ and $\Psi(v)=v_D$, with $C,D \subset \mathcal E_n$, we get
\begin{equation}\label{cambiospazio2}
\E(z_C v_D) = \sum_{A\subseteq \mathcal E_n}
P(A) \E_{A^c; \boldsymbol{\beta},\boldsymbol{0}}^{\mathrm{Is}}( z_C|_{z_i=0, \forall i \in A})
\E_{A; \boldsymbol{\alpha'},\boldsymbol{h'}}( v_D|_{v_{i}=\frac 12 , \forall i \in A^c})\,.
\end{equation}

The proof of the two inequalities \eqref{corr_in_1EG} and \eqref{corr_in_2EG} is an immediate application of the FKG inequality
relatively to $P$. 
Indeed, if $\boldsymbol{\alpha},\,\boldsymbol{h}\geq 0$,
then also $\boldsymbol{\b},\,\boldsymbol{\alpha'}, \boldsymbol{h'}\geq 0$,
and the conditions for the application of the FKG inequality are fulfilled:
\begin{itemize}
\item If $\boldsymbol{\b}\geq 0$, the GKS inequality for ferromagnetic Ising systems \cite{FV} guarantees that
the function $\E_{A^c; \boldsymbol{\beta},\boldsymbol{0}}^{\mathrm{Is}}( z_C|_{z_i=0, \forall i \in A})$
is non-increasing in $A$, for  any choice of $C \subseteq \mathcal E_n$.
%In particular, if  $A\subseteq B$,
%then $\E_{\mu_{A^c; \beta,0}}( z_C|_{z_i=0, \forall i \in A})
%\geq \E_{\mu_{B^c; \beta,0}}( z_C|_{z_i=0, \forall i \in B})$.
%
\item 
If $\boldsymbol{\alpha'},\,\boldsymbol{h'}\geq 0$ 
the function $\E_{A; \boldsymbol{\alpha'},\boldsymbol{h'}}( v_D|_{v_{i}=\frac 12 , \forall i \in A^c})$
is non-decreasing in $A$,  for any choice of  $D \subseteq \mathcal E_n$.
This is a consequence of the GKS inequality together with
Lemma \ref{lem_GKS-volume}.
Indeed, let $A\subseteq B$ and observe that
\be\label{passo0}
\begin{split}
&\E_{A; \boldsymbol{\alpha'},\boldsymbol{h'}}( v_D|_{v_{i}=\frac 12 , \forall i \in A^c})
= \frac{1}{2^{|D\cap A^c|}}
\E_{A; \boldsymbol{\alpha'},\boldsymbol{h'}}( v_{D\cap A})\\
&
\E_{B; \boldsymbol{\alpha'},\boldsymbol{h'}}( v_D|_{v_{i}=\frac 12 , \forall i \in B^c})
= \frac{1}{2^{|D\cap B^c|}}
\E_{B; \boldsymbol{\alpha'},\boldsymbol{h'}}( v_{D\cap B})\,.
\end{split}
\ee
Since $D\cap A\subseteq D\cap B$ by hypothesis, we can write
$v_{D\cap B}= v_{D\cap A}v_{D\cap (B\setminus A)}$ and hence,
from the GKS inequality and Lemma \ref{lem_GKS-volume},
\be\label{passo1}
\begin{split}
\E_{B; \boldsymbol{\alpha'},\boldsymbol{h'}}( v_{D\cap B})
&\geq \E_{B; \boldsymbol{\alpha'},\boldsymbol{h'}}( v_{D\cap A})
\E_{B; \boldsymbol{\alpha'},\boldsymbol{h'}}( v_{D\cap (B\setminus A)})\\
&\geq \E_{A; \boldsymbol{\alpha'},\boldsymbol{h'}}( v_{D\cap A})
\E_{A; \boldsymbol{\alpha'},\boldsymbol{h'}}( v_{D\cap (B\setminus A)})
\,.
\end{split}\ee
We now recall that for $\boldsymbol{\alpha'}\geq 0$ 
and $\boldsymbol{h'}\geq 0$, it holds that  $\E_{A; \boldsymbol{\alpha'},\boldsymbol{h'}}( v_i)\geq 1/2$, for all $i\in A$ and $A\subseteq \cE_n$.
Applying the GKS inequality to the second factor of the r.h.s of the above equation, 
and using this bound, we then get
\be\label{passo2}
\E_{A; \boldsymbol{\alpha'},\boldsymbol{h'}}(v_{D\cap (B\setminus A)})
\geq \prod_{i\in D\cap (B\setminus A)}
\E_{A; \boldsymbol{\alpha'},\boldsymbol{h'}}( v_i)
\geq  \frac{1}{2^{|D\cap(B\setminus A)|}}\,.
\ee
Putting together \eqref{passo0}-\eqref{passo2}, 
we conclude that  \be\label{GKS-wedge}
\E_{B; \boldsymbol{\alpha'},\boldsymbol{h'}}(v_D|_{v_{i}=\frac 12 , \forall i \in B^c})
\geq
\E_{A; \boldsymbol{\alpha'},\boldsymbol{h'}}( v_D|_{v_{i}=\frac 12 , \forall i \in A^c})\,.
\ee
\item If $\boldsymbol{\alpha} \geq 0$, the probability $P$, defined  in \eqref{misura_P} and acting on subsets of $\mathcal E_n$, 
satisfies the FKG lattice condition, namely
\be\label{FKG_P}
P(E)P(F)\leq P(E\cup F)P(E\cap F)\,,\qquad \forall E,F \subseteq \cE_n\,.
\ee
According to the definition of $P$, the inequality \eqref{FKG_P}
follows if the two inequalities
$$
Z_{E;\boldsymbol{\beta},\boldsymbol{0}}^{\mathrm{Is}}Z_{F;\boldsymbol{\beta},\boldsymbol{0}}^{\mathrm{Is}}
\leq
Z_{E\cup F;\boldsymbol{\beta},\boldsymbol{0}}^{\mathrm{Is}}Z_{E\cap F;\boldsymbol{\beta},\boldsymbol{0}}^{\mathrm{Is}}
$$
and
$$
Z_{E;\boldsymbol{\alpha'},\boldsymbol{h'}}Z_{F;\boldsymbol{\alpha'},\boldsymbol{h'}}
\leq
Z_{E\cup F;\boldsymbol{\alpha'},\boldsymbol{h'}}Z_{E\cap F;\boldsymbol{\alpha'},\boldsymbol{h'}}
$$
are simultaneously satisfied. The first inequality holds true  as a consequence  of the GKS inequality for Ising spin systems with $\boldsymbol{\beta} \geq 0$ and magnetic field $\boldsymbol{h} \geq 0$ (see \cite{L}, Lemma on p.~90).
The second inequality is instead verified thanks to Lemma~\ref{lem_zeta}.
\end{itemize}

Thus, as $P$ obeys the FKG lattice condition, and the functions $\E_{A^c; \boldsymbol{\beta},\boldsymbol{0}}^{\mathrm{Is}}( z_C|_{z_i=0, \forall i \in A})$ and $\E_{A; \boldsymbol{\alpha'},\boldsymbol{h'}}( v_D|_{v_{i}=\frac 12 , \forall i \in A^c})$
are respectively non-increasing and non-decreasing in $A$, 
from Eq. \eqref{cambiospazio2} we get
\begin{equation}\label{cambiospazio3}
\begin{split}
\E (z_C v_D )
&\leq \sum_{A\subseteq \mathcal E_n}
P(A) \E_{A^c; \boldsymbol{\beta},\boldsymbol{0}}^{\mathrm{Is}}( z_C|_{z_i=0, \forall i\in A})
\sum_{A\subseteq \mathcal E_n}
P(A)
\E_{A; \boldsymbol{\alpha'},\boldsymbol{h'}}( v_D|_{v_{i}=\frac 12 , \forall i \in A^c})\\
& = \E( z_C )\E( v_D )
 \,,
\end{split}
\end{equation}
providing inequality \eqref{corr_in_1EG}.
The inequality \eqref{corr_in_2EG} can be obtained in the same way by setting $\phi(z)=z_{C}z_{D}$,
so that  $g^{\psi}(A)=1$, and by observing that $f^{\phi}(A)$ is non-decreasing in $A$,
hence giving the reverse inequality.
\end{proof}

\begin{proof}[Proof of Theorem~\ref{Thm_GHS}.] The statement follows readily from Remark~\ref{chiusura} and Proposition~\ref{Prop:GHS-zv}.
\end{proof}

%%%%%%
\subsection{The GHS inequality for the two-star model} Let $\alpha, h \in \mathbb{R}$. Recall that the two-star model is obtained, as a particular case, by setting $\alpha_{ij} \equiv \alpha$, for all $i,j \in \mathcal{E}_n$, and $h_i \equiv h$, for all $i \in \mathcal{E}_n$, in the Hamiltonian \eqref{EW_H2}. This means that, whenever $\alpha, h \geq 0$, the GHS inequality holds true for the Gibbs measure $\mu_{n;\alpha,h}$, given in \eqref{eq:prob-exprg}. 

Observe that, by differentiating the free energy $f_{n; \a,h}$ w.r.t. $h$,
we get the following identities in terms of the Ursell functions 
(see also Remark \ref{rem-Ursell})
%\eqref{ursell}:
\[
\begin{split}
n^2 \partial_{h} f_{n;\a,h}&= \sum_{i\in\cE_n} u_1(i)\,,\quad
n^2 \partial_{h h} f_{n;\a,h}= \sum_{i,j\in\cE_n} u_2(i,j)\,,\quad\\
&n^2 \partial_{h h h} f_{n;\a,h}= \sum_{i,j,k\in\cE_n} u_3(i,j,k)\,,
\end{split}
\]
and so on. Thus, not only the sign of each Ursell function
provides a specific correlation inequality between
the random variables $x_i$'s, but also it  gives a definite
sign to a derivative of the free energy. 

A direct computation easily shows that, being a variance, the second order partial derivative of $f_{n;\alpha,h}$ w.r.t. $h$ is always non-negative. Thus, proving that \mbox{$u_2(i,j) \geq 0$} (GKS inequality) is useful to know the covariance between $x_i$ and $x_j$, but it is somehow irrelevant to the purpose of showing that the free energy is a convex function of~$h$. 
On the contrary, the GHS inequality ($u_3(i,j,k) \leq 0$) is of particular importance, as it implies  that the average edge density \eqref{averages} is a concave function of the parameter $h$
at any fixed size of the graph.

Explicitly, setting $m_n(\a,h):=\frac{\mathbb{E}_{n;\alpha,h}\left(E_n\right)}{n^2}$ and 
assuming that $\a, h\geq 0$, from the GKS and GHS inequalities we readily get 
\be\label{convex}
\partial_h m_{n} (\alpha, h) 
=
\partial_{h h} f_{n; \alpha, h} \geq 0
\,,\qquad
\partial_{h h} m_{n} (\alpha, h) 
=
\partial_{h h h} f_{n; \alpha, h} \leq 0\,.
\ee
Understanding the limiting behavior of the above derivatives has then a twofold  purpose.
On the one hand, it allows to infer properties  regarding the edge density and its limiting behavior;
for example, the existence of $\lim_{n\to+\infty}\partial_{h} m_{n} (\alpha, h)$ 
is fundamental for proving the standard central limit theorem in \eqref{CLT} (see \cite{BCM}).
On the other hand, it is crucial for detecting the occurence of phase transitions, 
that are generally associated with the emergence of  singularities in the infinite 
size free energy.
In particular, one can exploit convergence results on the derivatives of convex functions 
to guarantee that the limits and the derivatives w.r.t. the external field commute
(see \cite[Lemma V.7.5]{E}), and then obtain proper regularity conditions. 
Notice that this procedure can be seen  as an alternative approach 
to the investigation of the hypotheses that allow for the  application of the Lee-Yang theorem \cite{LY}.  However, in this respect, the convexity property \eqref{convex}
provides a more specific information that may enter in the characterization of further features of the model.

%==========================================================================
\section{Discussion on possible extensions}\label{Sect:gen}
%==========================================================================%
The results presented in Lemmas~\ref{lat_cond}--\ref{lem_zeta} can be extended to the general case where the Hamiltonian is a function of the homomorphism densities of an arbitrary collection of subgraphs of the graph $G$. In this Section we will briefly elaborate on this.

In the sequel, we will be dealing with the general Hamiltonian \eqref{Hamiltonian} and the corresponding Gibbs probability density \eqref{eq:prob-exprg}. We will denote by $\mathbb{E}_{n;\boldsymbol{\beta}}$ the relative expectation. Moreover, as a standard choice in the literature, we will set the subgraph $H_1$ to be an edge.

Going through the proof of Lemma \ref{lat_cond}, it is easy to understand
that the crucial condition for the validity of the FKG lattice condition
is inequality \eqref{cc-FKG}. When moving to the general setting we are adopting,
the analogous condition reads as
\begin{equation}\label{cc-FKGgen}
\exp \left\{ n^2 \sum_{j=2}^k \b_j t(H_j,G) \right\} \geq 1\,.
\end{equation}
As consequence, since the homomorphism densities are non-negative,
 the FKG lattice condition is in force whenever the parameters $\b_2,\ldots, \b_k$ are non-negative. Thus, we obtain the following result.
\begin{lemma}\label{lat_cond-g}
For all $\b_1\in\R$ and $\b_2,\ldots, \b_k\geq 0$, the Gibbs measure
$\mu_{n;\boldsymbol{\b}}$ fulfills the FKG lattice condition
\begin{equation}\label{FKG_ltg}
\mu_{n;\boldsymbol{\beta}} (x \vee y) \, \mu_{n;\boldsymbol{\beta}} (x \wedge y) 
\geq \mu_{n;\boldsymbol{\beta}} (x) \, \mu_{n;\boldsymbol{\beta}} (y)
\quad \text{ for } x,y \in \mathcal{A}_n.
\end{equation}
\end{lemma}

Two immediate consequences of Lemma \ref{lat_cond-g} are the positive correlation between increasing functions of the configuration and, in turn, the GKS inequality for the Gibbs measure $\mu_{n;\boldsymbol{\beta}}$. Specifically, for all $\beta_1 \in \mathbb{R}$ and 
$\beta_2, \dots, \beta_k \geq 0$, all increasing functions $f$ and $g$, 
and all $A,B\subseteq \cE_n$, it holds
\begin{align}
\E_{n;\boldsymbol{\beta}}(f g)
\geq \E_{n;\boldsymbol{\beta}}(f)\cdot
\E_{n;\boldsymbol{\beta}}(g) & &\mbox{(FKG inequality)}\label{FKG_g}\\[.3cm]
\E_{n;\boldsymbol{\beta}}(x_A x_B)
\geq \E_{n;\boldsymbol{\beta}}(x_A)\cdot
\E_{n;\boldsymbol{\beta}}(x_B), && \mbox{(GKS inequality)} \label{GKS_g}
\end{align}
where $x_C=\prod_{i \in C} x_i$, for $C\subseteq \mathcal{E}_n$.\\

An extension of Lemmas \ref{lem_GKS-volume} and \ref{lem_zeta}
 to this general context is also straightforward. 
However, while they were crucial to prove the GHS inequality for the two-star model,
they are quite irrelevant in the present setting, as the techniques used in 
Subsection \ref{subsec-GHS} can not be replicated.

%The main problem in extending previous argument, based on the duplication trick
%\eqref{ET_joint_measure}, comes from the breaking of identity \eqref{splitting}.
%Indeed, in this general setting, we can not reach factorized the joint measure of the
%dublicate 

Indeed, when dealing with a generic exponential random graph, the trick of variable duplication 
\eqref{ET_joint_measure} does not work. The problem is twofold. On the one hand, in general the decomposition \eqref{splitting} fails to exist, as mixed terms remain. Thus, it is not possible to factorize the joint measure of the doubled model and then characterize correlations exploiting averages over an Ising and an ERG subsystem (see \eqref{mediafunzioni}). On the other hand, even if the joint measure factored out and the analog of \eqref{mediafunzioni} were available, to conclude we would need FKG and GKS inequalities for Ising models with multi-body interactions, that are not known. 

However, if we specify the Ursell function $u_3(i,j,k)$, given in \eqref{u3},
 in the special cases when $i=j=k$ and $i=j\neq k$, 
we obtain respectively 
\[
 \mathbb{E}_{n;\boldsymbol{\beta}}(x_i) \left(1-\mathbb{E}_{n;\boldsymbol{\beta}}(x_i)\right) (1-2\mathbb{E}_{n;\boldsymbol{\beta}}(x_i))
\]
and 
\[
 \text{Cov}_{n;\boldsymbol{\beta}}(x_i,x_j) (1-2\mathbb{E}_{n;\boldsymbol{\beta}}(x_i)) .
\]
%However, if we specify the Ursell function $u_3(i,j,k)$ in the special cases when $i=j=k$ and $i=j\neq k$, we obtain respectively 

%\[
%\sum_{i \in \mathcal{E}_n} \mathbb{E}_{n;\boldsymbol{\beta}}(x_i) \left(1-\mathbb{E}_{n;\boldsymbol{\beta}}(x_i)\right) (1-2\mathbb{E}_{n;\boldsymbol{\beta}}(x_i))
%\]
%and 
%\[
%\sum_{i,j \in \mathcal{E}_n} \text{Cov}_{n;\boldsymbol{\beta}}(x_i,x_j) (1-2\mathbb{E}_{n;\boldsymbol{\beta}}(x_i)) .
%\]
Since $\text{Cov}_{n;\boldsymbol{\beta}}(x_i,x_j) \geq 0$, due to the GKS inequality \eqref{GKS_g}, we can conjecture that the necessary and sufficient condition for the GHS inequality to hold in the present general setting is again only the requirement $\mathbb{E}_{n;\boldsymbol{\beta}}(x_i) \geq 1/2$.

\end{document}